\documentclass[a4paper,11pt]{amsart}
\usepackage{graphicx}
\usepackage[ansinew]{inputenc}    % accents i altres
\usepackage[all]{xy}

\usepackage{xcolor}%para el color de letras
\definecolor{colordelink}{rgb}{0,0,0.50}
\definecolor{colordecite}{rgb}{0,0.5,0}
\definecolor{colordeurl}{rgb}{0,0.41,0.5}
\usepackage[pagebackref=true]{hyperref}% Para linkear referencias. COMPILAR EN PDF!!! Si tal quitar colorlinks https://tex.stackexchange.com/questions/50747/options-for-appearance-of-links-in-hyperref 
\hypersetup{
    colorlinks=true,
    linkcolor=colordelink,
    citecolor=colordecite,
		urlcolor=colordeurl,
}

\usepackage[capitalize,nameinlink,noabbrev]{cleveref}% to emulate \autoref style

\usepackage{todonotes}

\usepackage{amssymb,amscd}

\usepackage{mathrsfs,mathtools}
\usepackage[shortlabels]{enumitem}%Para las opciones de enumerate como [(i)], shortlabels es para escribir i en lugar de roman etc.
\usepackage{tikz-cd}%Da problemas con graphicx Para diagramas conmutatvos y ese tipo de cosas: http://ctan.math.washington.edu/tex-archive/graphics/pgf/contrib/tikz-cd/tikz-cd-doc.pdf

\theoremstyle{definition}

\theoremstyle{definition}

\theoremstyle{definition}

\def\C{\mathbb C}

\def\dim{\operatorname{dim}}

\def\codim{\operatorname{codim}}

\def\id{\operatorname{id}}

\def\im{\operatorname{Im}}

\def\Iso{\operatorname{Iso}}
\def\Der{\operatorname{Der}}

 \newcommand{\medpar}[1]{\big(#1\big)}

\newcommand{\CC}{\mathbb{C}}

\newcommand{\ZZ}{\mathbb{Z}}
\newcommand{\eqA}{\mathscr{A}}

\def\aecodim{\eqA_e\text{-}\codim}

\newcommand{\GS}[2]{(\mathbb{C}^{#1},S)\rightarrow(\mathbb{C}^{#2},0)}

\theoremstyle{plain}
\newtheorem{theorem}{Theorem}[section]
\newtheorem{lemma}[theorem]{Lemma}
\newtheorem{corollary}[theorem]{Corollary}
\newtheorem{proposition}[theorem]{Proposition}

\theoremstyle{definition}
\newtheorem{definition}[theorem]{Definition}
\newtheorem{conjecture}[theorem]{Conjecture}
\newtheorem{example}[theorem]{Example}

\newtheorem{remark}[theorem]{Remark}

\begin{document}

\author{R. Giménez Conejero and
J.J.~Nu\~no-Ballesteros}

\title[A weak version of Mond's conjecture]
{A weak version of  Mond's conjecture}

\address{%Departament de Matemàtiques,
%Universitat de Val\`encia, Campus de Burjassot, 46100 Burjassot
%SPAIN
Alfr\'ed R\'enyi Institute of Mathematics, Re\'altanoda utca 13-15,
H-1053 Budapest, 
Hungary
}
\email{Roberto.Gimenez@uv.es}
\address{Departament de Matem\`atiques,
Universitat de Val\`encia, Campus de Burjassot, 46100 Burjassot
SPAIN. \newline Departamento de Matemática, Universidade Federal da Paraíba
		CEP 58051-900, João Pessoa - PB, BRAZIL}
\email{Juan.Nuno@uv.es}

\thanks{Grant PGC2018-094889-B-100 funded by MCIN/AEI/ 10.13039/501100011033 and by ``ERDF A way of making Europe''.}


\begin{abstract} 
We prove that a map germ $f:(\CC^n,S)\to(\CC^{n+1},0)$ with isolated instability is stable if and only if $\mu_I(f)=0$, where $\mu_I(f)$ is the image Milnor number defined by Mond. In a previous paper we proved this result with the additional assumption that $f$ has corank one. The proof here is also valid for corank $\ge 2$, provided that $(n,n+1)$ are nice dimensions in Mather's sense (so $\mu_I(f)$ is well defined). Our result can be seen as a weak version of a conjecture by Mond, which says that the $\eqA_e$-codimension of $f$ is $\le \mu_I(f)$, with equality if $f$ is weighted homogeneous. As an application, we deduce that the bifurcation set of a versal unfolding of $f$ is a hypersurface.
\end{abstract}

\maketitle
\bigskip
\section{Introduction}

In the context of Thom-Mather theory (deformation theory) for map germs $f:(\CC^n,S)\to(\CC^{p},0)$, the cases $p<n+1$, $p=n+1$ and $p>n+1$ present different traits and behaviours. In particular, the case $p=n+1$ has several open questions that are understood for the other cases. Here we solve two of them: one regarding the homotopy type of a \textit{generic fiber} and the other regarding the \textit{bifurcation set}.
\newline

A hypersurface $(X,0)$ with isolated singularity has a well known invariant given by its Milnor fiber, the Milnor number $\mu(X,0)$. D. Mond introduced a similar invariant for map germs $f:(\CC^n,S)\to(\CC^{n+1},0)$ with isolated instability in \cite{Mond1991}. Indeed, the image of a stable perturbation of $f$ has the homotopy type of a wedge of spheres of dimension $n$, whose number is independent of the stable perturbation. This number of spheres is the image Milnor number of $f$, denoted by $\mu_I(f)$, and its stable perturbation plays the role of Milnor fiber of a hypersurface as before. However, this is only well defined when $f$ has corank one or, alternatively, $(n,n+1)$ are nice dimensions in Mather's sense (cf. \cite{Mather1971}).

On the other hand, the $\eqA_e$-codimension of a germ as above is the equivalent of the Tjurina number $\tau(X)$ of a hypersurface with isolated singularity, because they control the space of perturbations of map germs and germs of hypersurfaces, respectively.

The easy relation between the Tjurina and Milnor number for hypersurfaces inspired D. Mond to conjecture in \cite{Mond1991} that the relation is reproduced in the case of map germs with isolated instability. More precisely:

\begin{conjecture}[Mond's conjecture]\label{MondConjecture}
Given a germ $f:(\CC^n,S)\to(\CC^{n+1},0)$ with finite $\eqA_e$-codimension such that it has corank one or $(n,n+1)$ are nice dimensions (i.e., $n<15$),
$$ \mu_I(f)\geq\aecodim(f), $$
with equality in the quasi-homogeneous case.
\end{conjecture}

This question remains open in general (see {\cite[Theorem 4.2]{deJong1991} and \cite[Theorem 2.3]{Mond1995}} for the cases $n=1,2$).

In this paper we show a basic result on both objects of map germs. In \cref{Sec: mui} we deal with the image Milnor number an prove that it is positive if the germ is, indeed, unstable. This is a weak version of \cref{MondConjecture} (i.e., when $\mu_I(f)=0$ we have $\aecodim(f)=0$) but it implies this conjecture in the case of $\eqA_e$-codimension one and the case of \textit{augmentations} of such germs. In order to prove these results we also introduce new objects such as $LC(G)$, a version of Saito's characteristic variety $LC(\mathcal{X})$, to control the image Milnor number.
Furthermore, a deep understanding of $LC(G)$ allows us to prove that, if it is Cohen-Macaulay, Mond's conjecture is true for germs with a one parameter stable unfolding.
\newline

The bifurcation set is the set of parameters of a versal unfolding such that the corresponding perturbations have instabilities. In \cref{sec: bif} we prove that the bifurcation set of a germ $f:(\CC^n,S)\to(\CC^{n+1},0)$ is a pure dimensional hypersurface. Its relevance is also put into perspective by explaining that the methods used to prove a similar fact in other settings fail in our case. Furthermore, a profound understanding of the space of unstable perturbations could lead to a proof of Mond's conjecture (see the program to solve this conjecture in the first author's thesis, \cite[Section 7.3]{Robertothesis}). Knowing that the bifurcation set is of codimension one is the first step in that direction.

We refer to the modern reference \cite{Mond2020} for the definitions and properties about singularities of mappings such as stability, finite determinacy, versal unfoldings, and other basic concepts.

\section{Unstable germs and the image Milnor number}\label{Sec: mui}
Throughout this section we assume that $f\colon(\CC^{n},S)\to(\CC^{n+1},0)$ is an $\eqA$-finite germ such that it has corank one or $(n,n+1)$ are nice dimensions.

%We denote by $(X,0)$ the image of $f$ in $(\CC^{n+1},0)$. 
We fix a representative $f\colon U\to V$ which is a finite mapping (i.e., finite-to-one and closed) and denote by $X=f(U)$ its image. The germ of $f$ at the point $y\in X$ is denoted by
$$(f)_{y}:\medpar{\CC^n,f^{-1}(y)}\to(\CC^{n+1},y).$$

There is a natural stratification on an $\eqA$-finite germ $f$, which is straightforward in the case of a stable germ.

\begin{definition}\label{stratification_stable_types}
Suppose $f$ is stable. The \textit{isosingular locus} of $f$ at a point $y_0$ in the image $X$ is the set
$$  \Iso\big(f;y_0\big)\coloneqq\left\{y\in X: (f)_y \text{ is }\eqA\text{-equivalent to } (f)_{y_0}\right\}.$$ Moreover, the \textit{stratification by stable types} of the image of $f$ is given by the isosingular loci of all points $y_0\in X$. There is an induced stratification in the source of $f$, these two stratifications are the \textit{stratification by stable types of $f$}.
\end{definition}

\begin{figure}
	\centering
		\includegraphics[width=1.00\textwidth]{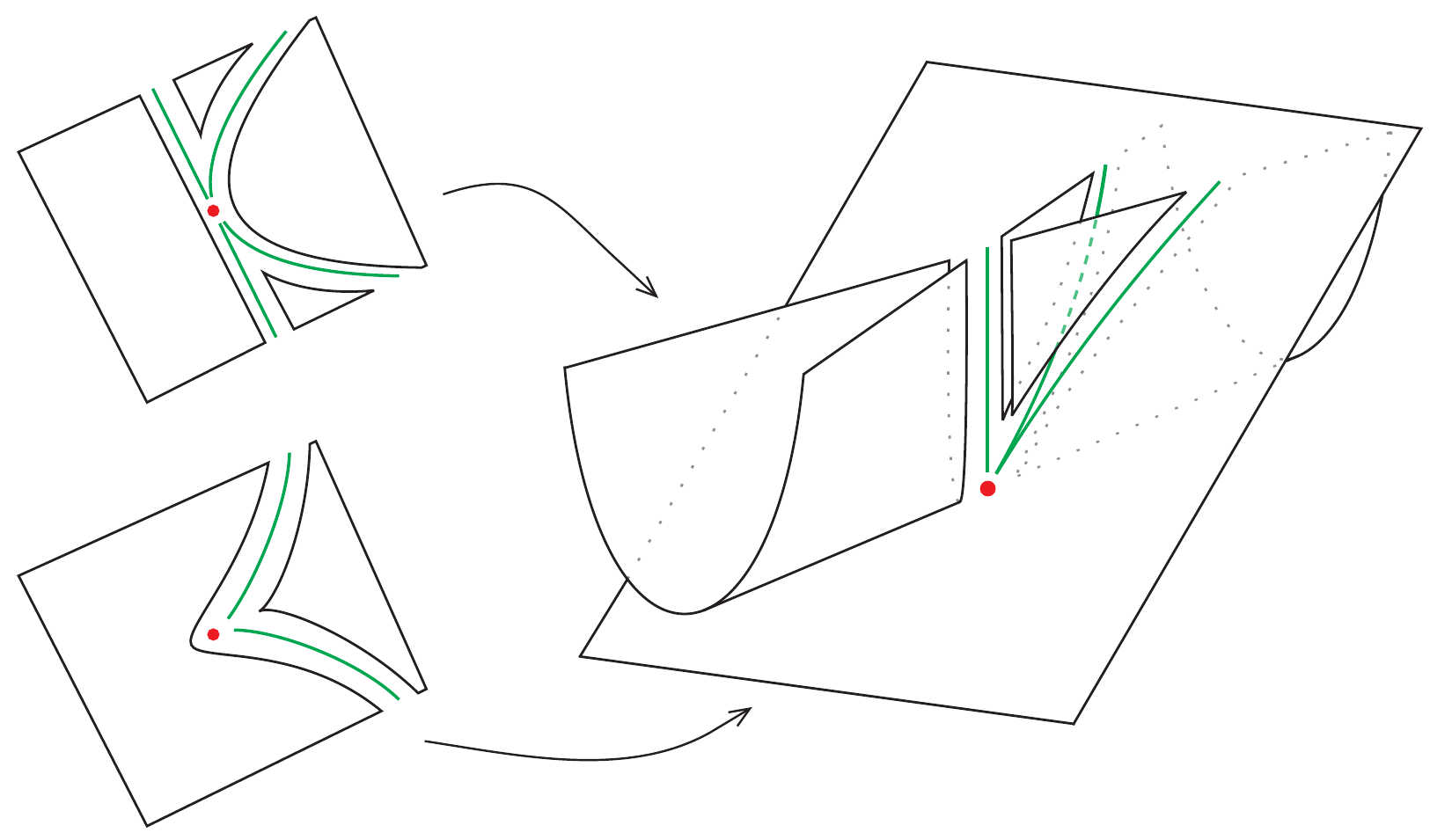}
	\caption{Stratification by stable types of an unstable bigerm given by a crosscap and an immersion. Observe that the transverse double points are in the same stratum.}
	\label{fig:StratificationStableTypesDef}
\end{figure}

We refer to \cite[Definition 7.2]{Mond2020} for details and properties of this stratification. In fact, the hypothesis that $f$ has corank one or
$(n,n+1)$ are in the nice dimensions guarantees that we have a finite number of strata.
The name of \textsl{stable types} is taken because, obviously, this stratification identifies stable singularities of the same $\eqA$-class (cf. \cref{fig:StratificationStableTypesDef}). Indeed, there is an analogous definition of the \textit{stratification by stable types} of a locally stable map.

In contrast, it could happen that $f:\GS{n}{n+1}$ has instabilities. However, it has isolated instability if it is $\eqA$-finite (by Mather-Gaffney criterion). We can 
assume $f\colon U-S\to V-0$ is locally stable and extend its stratification by stable types just by adding the strata $S$ and $0$ in source and target, respectively. We call this stratification \textit{stratification by stable types of $f$} as well (see \cref{fig:StratificationStableTypesDef}).
\newline

Let $F(x,s)=(f_s(x),s)$ be a stabilisation of $f\colon(\CC^{n},S)\to(\CC^{n+1},0)$. Then, $F$ is also $\eqA$-finite, considered as a germ $(\CC^n\times\CC,S\times 0)\to(\CC^{n+1}\times\CC,0)$. We denote by $(\mathcal X,0)$ the image of $F$ in $(\CC^{n+1}\times\CC,0)$ and by $\pi\colon(\mathcal X,0)\to(\CC,0)$ the restriction of the projection onto the parameter space. For each $s$ close to $0$ in $\CC$, the fibre $\pi^{-1}(s)$ is $X_s$, the image of $f_s$.

\begin{figure}
	\centering
		\includegraphics[width=1.00\textwidth]{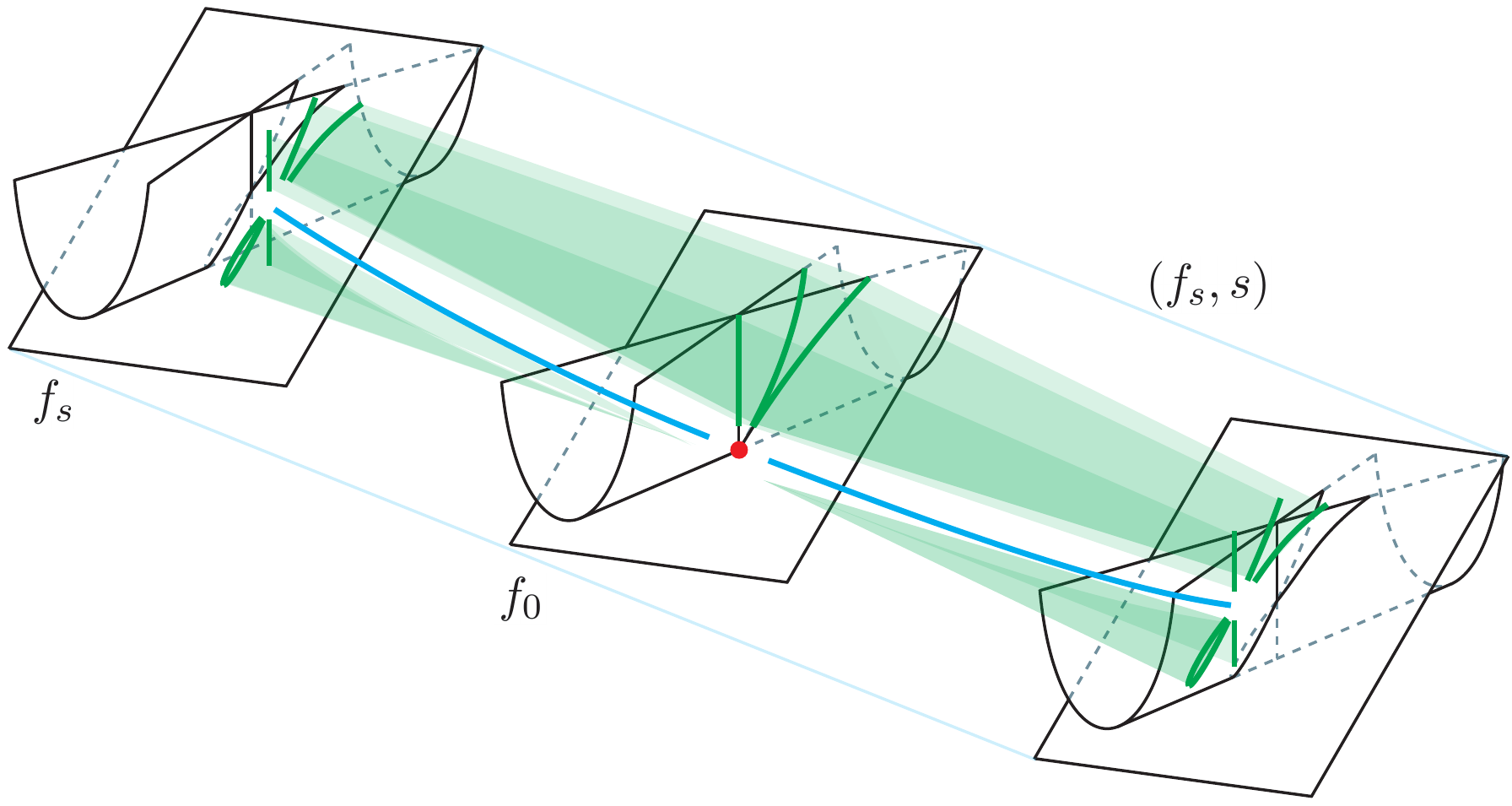}
	\caption{Stratification by stable types of the image of an unfolding $(f_s,s)$ of the unstable bigerm $f_0$ given by a crosscap and an immersion. Observe that, in contrast to \cref{fig:StratificationStableTypesDef}, the stratum of dimension zero in now stable as a germ from $\CC^3$ to $\CC^4$. In a general unfolding this stratum could be unstable.}
	\label{fig:StratificationUnfoldingDef}
\end{figure}

We remark that $F$ still has a well defined stratification by stable types on its image $\mathcal X$, although now $(n+1,n+2)$ could be in the boundary of the nice dimensions. In fact, if $(y,s)\in \mathcal X-0$, then the germ of $f_s$ at $y$ is stable. Therefore, the germ of $F$ at $(y,s)$ is a trivial unfolding of the germ of $f_s$ at $y$. Since $(n,n+1)$ are nice dimensions or the germ has corank one, only a finite number of stable types of germs from $\CC^n$ to $\CC^{n+1}$ can appear. Thus, the family of isosingular loci $\Iso\big(F;(y,s)\big)$ at points $(y,s)\in\mathcal X-\{0\}$ is finite. By adding the origin if necessary as a new stratum, we get the stratification by stable types of $\mathcal X$ (see \cref{fig:StratificationUnfoldingDef}).

%Furthermore, the germ of $f_s$ at the point $y\in X_s$ is denoted as
%$$(f_s)_{y}:\medpar{\CC^n,f_s^{-1}(y)}\to(\CC^{n+1},y),$$
%%where $\Sigma(f)$ is the locus of critical points of $f$. 
%\textcolor{red}{aquí he quitado $\Sigma(f)$}

\begin{lemma} The (connected components of the) stratification by stable types of the hypersurface $(\mathcal X,0)$ coincides with the logarithmic stratification. In particular, $(\mathcal X,0)$ is holonomic.
\end{lemma}

\begin{proof} Given two points $(y,s)$ and $(y',s')$ in $\mathcal X$ in the same connected component of the stratification by stable types, the germs of $F$ at $(y,s)$ and $(y',s')$ are $\eqA$-equivalent. In particular, the germs of $\mathcal X$ at $(y,s)$ and $(y',s')$ are equivalent by a biholomorphism $\Psi$ of the ambient space $\CC^{n+1}\times\CC$. Thus, by connectivity $(y,s)$ and $(y',s')$ belong to the same stratum in the logarithmic stratification of $\mathcal X$.

But the converse also holds. Indeed, the germ $F$ at each point is the normalisation of $\mathcal X$ at that point. 
If the germs of $\mathcal X$ at $(y,s)$ and $(y',s')$ are equivalent by a biholomorphism $\Psi$ of $\CC^{n+1}\times\CC$, there exists a unique biholomorphism $\Phi$ in $\CC^n\times\CC$ so that $(\Phi,\Psi)$ gives an $\eqA$-equivalence between the germs of $F$ at $(y,s)$ and $(y',s')$. Modulo connectivity, this shows that the stratification by stable types coincides with the logarithmic stratification.%That strat is locally a product, that is why this works
\end{proof}

\begin{lemma}\label{isolated} The projection $\pi\colon(\mathcal X,0)\to(\CC,0)$ has isolated critical points in the stratified sense. Moreover, it has a critical point if, and only if, $f$ is not stable.
\end{lemma}

\begin{proof} As above, we fix a representative $\mathcal X$ and take $(y,s)\in\mathcal X-\{0\}$. Since the germ of $f_s$ at $y$ is stable, the germ of $F$ at $(y,s)$ is a trivial unfolding of the germ of $f_s$ at $y$. Hence, there exist biholomorphisms $\Phi$ and $\Psi$ in $\CC^n\times\C$ and $\CC^{n+1}\times\C$, respectively, which are unfoldings of the identities and such that $\Psi\circ F\circ\Phi^{-1}=f_s\times\id_{\CC}$ in a neighbourhood of $(y,s)$.
This gives a commutative diagram
$$\begin{tikzcd}
\big(\mathcal X,(y,s)\big) \arrow[r,"\pi"] \ar[d,"\Psi"',"\sim"  {anchor= north, rotate=90, inner sep=.6mm}]&(\CC,s)\\
\big(X_s\times\CC,(y,s)\big)\ar[ru,"\pi_2"'] &
\end{tikzcd},$$
so $(y,s)$ is a regular point of $\pi$ in the stratifed sense.

Assume now that the origin $0$ is also a regular point of $\pi$ in the stratifed sense. Let $S_0=\Iso(F;0)$ be the stratum of $\mathcal X$ which contains $0$. Obviously, we must have $S_0\ne\{0\}$, and since $F$ is stable at any point $(y,s)\in S_0-\{0\}$, $F$ is also stable at $0$. As $0\in S_0$ is a regular point of the restriction $\pi\colon S_0\to\C$, we deduce that %$0\in\C$ is a regular value of $\pi\colon S\to\C$. This implies that
the hyperplane $\pi^{-1}(0)=\CC^{n+1}\times\{0\}$ is transverse to $S_0$. By \cite[Proposition 2.22]{NunoBallesteros2018}, $f_0=f$ is stable at $0$. The converse is obvious, as any unfolding is trivial and has the form $F(x,t)=\big(f_t(x),t\big)$ up to $\eqA$-equivalence.
\end{proof}

\begin{example}
In \cref{fig:StratificationUnfoldingDef} we have represented an unfolding $(f_s,s)$ of the unstable bigerm $f_0:\big(\CC^2,\left\{p,q\right\}\big)\to(\CC^3,0)$ given by a crosscap and an immersion (see also \cref{fig:StratificationStableTypesDef}), where $f_s$ is a stable perturbation of $f_0$ (i.e., $(f_s,s)$ is a stabilisation). 

It is easy to see that \cref{isolated} holds in this case. Indeed, the unique stratum of dimension zero is also the unique critical point of the projection to the parameter $s$.

The fact that this point is stable as a germ induced by the unfolding $(f_s,s)$ (i.e., as a germ from $\CC^3$ to $\CC^4$) plays no role in this fact. In general, if we begin with an unstable $f_0$ then, in the unfolding, we need to add as a stratum the point where the instability of $f_0$ is located, regardless whether it is stable or unstable as a germ induced by the unfolding.
\end{example}

Let $\theta_{n+2}$ be the $\mathscr O_{n+2}$-module of germs of vector fields on $\CC^{n+2}\equiv \CC^{n+1}\times\CC$ at the origin. We denote by $\Der(-\log\mathcal X)$ the submodule of \textit{logarithmic vector fields}. We recall that $\xi\in\Der(-\log\mathcal X)$ if and only if $\xi_p\in T_p\mathcal X$, for all $p\in\mathcal X_{reg}$, the regular part of $\mathcal X$. Equivalently, $\xi\in\Der(-\log\mathcal X)$ if and only if $dG(\xi)\in(G)$, where $G\in\mathscr O_{n+2}$ such that $G=0$ is a reduced equation of $\mathcal X$.

Take a representative $\mathcal X$ in some open neighbourhood $U$ of the origin in $\CC^{n+2}$. We extend the stratification of $\mathcal X$ to $U$ by adding the open stratum $U-\mathcal X$.
The projection $\pi\colon U\to \CC$ has also an isolated stratified critical point at the origin and, hence, the \textit{Bruce-Roberts number}
\[
\mu_{BR}(\mathcal X,\pi):=\dim_\CC\frac{\mathscr O_{n+2}}{d\pi\big(\Der(-\log\mathcal X)\big)}
\]
is always finite (see \cite[Definition 2.4]{Bruce1988} and the previous comments) and is not zero when $\pi$ has a critical point at the origin, that is, when $f$ is not stable (by \cref{isolated}). It seems natural to ask about the relationship between this number and $\mu_I(f)$, which gives the number of vanishing cycles of the fiber $X_s=\pi^{-1}(s)$. However, the following example shows that these two numbers are not equal in general.

\begin{example}\label{exam1} Let $f\colon(\C^2,0)\to(\C^3,0)$ be given by $f(x,y)=(x^2,y^2,xy+x^3+y^3)$ with stabilisation
\[
F(x,y,s)=(x^2+sy,y^2-sx,xy+y^3+x^3+s(x-y),s).
\]
We have $\mu_I(f)=7$ (see \cite[Section 3.1]{Marar2008}) and a computation with \textsc{Singular} (see \cite{DGPS}) gives that $\mu_{BR}(\mathcal X,\pi)=8$.
\end{example}

Instead of $\Der(-\log\mathcal X)$ we will consider the submodule $\Der(-\log G)$, defined as the set of vector fields $\xi$ such that $dG(\xi)=0$. Roughly speaking, the difference between these modules is disregarding the fiber $G^{-1}(0)=\mathcal X$ as special and considering the tangency at every fiber of $G$. Obviously, we have the inclusion $\Der(-\log G)\subset \Der(-\log\mathcal X)$. Moreover, when $\mathcal X$ is weighted homogeneous, we also have that
\begin{equation}\label{euler}
\Der(-\log\mathcal X)=\Der(-\log G)\oplus \mathscr O_{n+2}\cdot\{\epsilon\},
\end{equation}
where $\epsilon$ is the Euler vector field
$$\epsilon= w_1 \frac{\partial}{\partial x_1}+\cdots+w_{n+1}\frac{\partial}{\partial x_{n+1}}+w_s\frac{\partial}{\partial s},$$
$w_\bullet$ denoting the corresponding weights of the variables. Moreover, since $d\pi(\epsilon)=w_{s}s$, we obtain the equality
\begin{equation}\label{euler2}
d\pi\big(\Der(-\log\mathcal X)\big)=d\pi\big(\Der(-\log G)\big)+(s).
\end{equation}
\medskip

\begin{figure}
	\centering
		\includegraphics[width=1.00\textwidth]{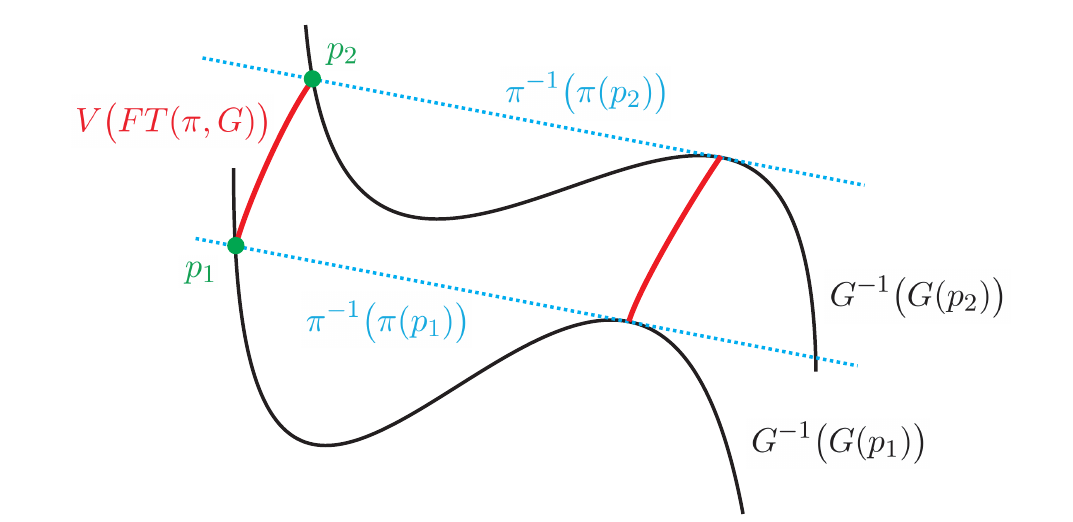}
	\caption{Representation of the set of zeros of the ideal $ FT(\pi,G)= d\pi\big(\Der(-\log G)\big)$. Observe that it has two components in this representation.}
	\label{fig:FT}
\end{figure}

The ideal $d\pi\big(\Der(-\log G)\big)\subset \mathscr{O}_{n+2}$ contains relevant information we want to analyse. Observe that if, at a point $p$, the tangent space of the fiber of $G(p)$ coincides with the tangent space of the fiber of $\pi(p)$, i.e., if $$ T_p G^{-1}\big(G(p)\big)=T_p \pi^{-1}\big(\pi(p)\big),$$ then $p$ will be in the set of zeros of $d\pi\left(\Der(-\log G)\right)\subset \mathscr{O}_{n+2}$ (see \cref{fig:FT}). Hence, at least,
 this ideal contains the information of
 the \textsl{failure of the transversality} between the fibers of $\pi$ and $G$. For this reason, we use the notation $$ FT(\pi,G)\coloneqq d\pi\big(\Der(-\log G)\big).$$ Note, however, that if we consider the points where the fibers of $G$ are not smooth this does not work, so $FT(\pi,G)$ refines this notion.

%$$ \textcolor[rgb]{0,0.53,0.22}{p_1\ p_2}\ \textcolor[rgb]{0,0.62,0.76}{\pi^{-1}\big(\pi(p_1)\big)\ \pi^{-1}\big(\pi(p_2)\big)}\ G^{-1}\big(G(p_1)\big)\ G^{-1}\big(G(p_2)\big)\ \textcolor[rgb]{1,0,0}{V\big(FT(\pi,G)\big)}  $$

\medskip

In the following lemma, we consider the general case where $\mathcal X$ is not necessarily weighted homogeneous. 
For each $s\in\CC$ close enough to the origin, we denote by $g_s$ the function $g_s(y)=G(y,s)$.

\begin{lemma} \label{V(I)}
The zero locus $V\big(FT(\pi,G)\big)$ of $FT(\pi,G)$ in $(\CC^{n+2},0)$ is the set-germ of points $(y,s)$ such that either:
\begin{enumerate}
\item $(y,s)=0$ and $f$ is not stable, or 
\item $G(y,s)\ne0$ and $y$ is a critical point of $g_s$. 
\end{enumerate}
Moreover, in the second case, $FT(\pi,G)$ is generated by $\frac{\partial G}{\partial y_1},\dots,\frac{\partial G}{\partial y_{n+1}}$ in a neighbourhood of $(y,s)$.
\end{lemma}

\begin{proof} By the curve selection lemma, $G$ has isolated critical value at the origin. We fix a representative in some open neighbourhood $U$ of the origin in $\CC^{n+2}$ such that $f_s$ is stable at $y$, for all $(y,s)\in \mathcal X-\{0\}$ and $0$ is the only critical value of $G$.

%Let $p=(y_0,s_0)\in U$. We first consider the case $G(p)=0$, so $p\in\mathcal X$. If $p\ne0$, then $p$ is not a critical point of $\pi\colon U\to\CC$ (in the stratified sense) and hence, $p\notin V(d\pi\left(\Der(-\log\mathcal X)\right))$. 
%
%On the other hand, $f_{s_0}$ is stable at $y_0$ so the germ of $F$ at $p$ is a trivial unfolding of the germ of $f_{s_0}$ at $y_0$. Since $(n,n+1)$ are nice dimensions, the germ of $f_{s_0}$ at $y_0$ is weighted homogeneous, up to $\eqA$-equivalence (see \cite[Theorem 7.6]{Mond2020}). In particular, $\mathcal X$ is weighted homogeneous in a neighbourhood of $p$, up to a coordinate change which preserves the parameter $s$. It follows from \eqref{euler2} that 
%\[
%d\pi\left(\Der(-\log\mathcal X)\right)=d\pi\left(\Der(-\log G)\right)+(s-s_0)=I+(s-s_0),
%\]
%in a neighbourhood of $p$. But this implies that $p\notin V(I)$.
%
%When $p=0$ and $f$ is stable at $0$ then we use the same argument as above and arrive to that $0\notin V(I)$. If $p=0$ and $f$ is not stable at $0$, then 
%$0$ is a critical point of $\pi\colon U\to\CC$ (in the stratified sense) and hence, 
%\[
%0\in V(d\pi\left(\Der(-\log\mathcal X)\right)\subset V(I).
%\]
Let $p=(y_0,s_0)\in U$. We first consider the case $G(p)=0$, so $p\in\mathcal X$. 

If $f_{s_0}$ is stable at $y_0$ (in particular, if $p\neq0$), then $p$ is not a critical point of $\pi\colon U\to\CC$, by \cref{isolated}. Hence, $p\notin V\big(d\pi\big(\Der(-\log\mathcal X)\big)\big)$. On the other hand, at $p$, $F$ is a trivial unfolding of the germ $(f_{s_0})_{y_0}$, as it is stable. Since we deal with corank one germs or $(n,n+1)$ are nice dimensions, the germ $(f_{s_0})_{y_0}$ is weighted homogeneous, up to $\eqA$-equivalence (see \cite[Theorem 7.6]{Mond2020}). In particular, $\mathcal X$ is weighted homogeneous in a neighbourhood of $p$, up to a coordinate change which preserves the parameter $s$. It follows from \cref{euler2} that 
\[
d\pi\big(\Der(-\log\mathcal X)\big)=d\pi\big(\Der(-\log G)\big)+(s-s_0)=FT(\pi,G)+(s-s_0),
\]
in a neighbourhood of $p$. Hence, $p\notin V\big(FT(\pi,G)\big)$ because $p\in V(s-s_0)$.

If $p=0$ and $f$ is unstable at that point, then 
$0$ is a critical point of $\pi\colon U\to\CC$ (by \cref{isolated}) and, hence, 
\[
0\in V(d\pi\big(\Der(-\log\mathcal X)\big)\subset V\big(FT(\pi,G)\big).
\]

\medskip

Now we consider the case $G(p)\ne0$. By assumption, $G$ is regular at $p$. Assume that $\frac{\partial G}{\partial y_i}(p)\ne0$, for some $i=1,\dots,n+1$. We can suppose, for simplicity, that $i=1$. The map 
\[
\Phi(y,s)=\big(G(y,s),y_2,\dots,y_{n+1},s\big)
\] 
is a biholomorphism in a neighbourhood of $p$ and 
\[
G\circ\Phi^{-1}(z,y_2\dots,y_{n+1},s)=z,
\]
in a neighbourhood of $\Phi(p)$. The module $\Der\big(-\log (G\circ\Phi^{-1})\big)$ is generated at $\Phi(p)$ by the vector fields $\frac{\partial}{\partial y_2},\dots,\frac{\partial}{\partial y_{n+1}},\frac{\partial}{\partial s}$. Composing by the differential of $\Phi$ we reverse the coordinate change, so $\Der(-\log G)$ is generated at $p$ by
\begin{align*}
\hspace{2cm}\xi_i&=\frac{\partial}{\partial y_{i}}+\frac{\partial G}{\partial y_i}\frac{\partial}{\partial y_{1}},\quad i=2,\dots,n+1,\\
\xi_{n+2}&=\frac{\partial}{\partial s}+\frac{\partial G}{\partial s}\frac{\partial}{\partial y_{1}}.
\end{align*}
We have $d\pi(\xi_{n+2})=1$, so $p\notin V\big(FT(\pi,G)\big)$. 

The case $\frac{\partial G}{\partial s}(p)\ne0$ has to be analysed separately. We proceed analogously and arrive to that $\Der(-\log G)$ is generated at $p$ by 
\[
\eta_i=\frac{\partial}{\partial y_{i}}+\frac{\partial G}{\partial y_i}\frac{\partial}{\partial s},\quad i=1,\dots,n+1.
\]
In this case $d\pi(\eta_i)=\frac{\partial G}{\partial y_i}$. Hence, $FT(\pi,G)$ is generated at $p$ by $\frac{\partial G}{\partial y_1},\dots,\frac{\partial G}{\partial y_{n+1}}$, and $p\in V\big(FT(\pi,G)\big)$ if and only if $y_0$ is a critical point of $g_{s_0}$.
\end{proof}

\begin{corollary}\label{finite} The number
\[
\dim_\CC\frac{\mathscr O_{n+2}}{FT(\pi,G)+(s)}
\]
is always finite and is not zero if and only if $f$ is not stable. 
\end{corollary}

\begin{proof} By shrinking the neighbourhood $U$ if necessary, we can assume that $0$ is the only critical value of $g_0$. Hence, by \cref{V(I)}, we have that 
\[
V\big(FT(\pi,G)+(s)\big)\subset\{0\},
\]
with equality if and only if $f$ is not stable. The results follows now from the analytic Nullstellensatz.
\end{proof}

\begin{example} 
In \cref{exam1}, we have 
\[
\dim_\CC\frac{\mathscr O_{n+2}}{FT(\pi,G)+(s)}=7,
\]
which coincides with $\mu_I(f)$.

\end{example}

 In the next theorem, we consider $\mathscr O_{n+2}/FT(\pi,G)$ as an $\mathscr O_1$-module via the morphism $\pi^*\colon\mathscr O_1\to\mathscr O_{n+2}$. By \cref{finite}, $\mathscr O_{n+2}/FT(\pi,G)$ is always finitely generated over $\mathscr O_1$.
%Let $I=d\pi\left(\Der(-\log G)\right)$. In the next theorem, we consider $\mathscr O_{n+2}/I$ as an $\mathscr O_1$-module via the morphism $\pi^*\colon\mathscr O_1\to\mathscr O_{n+2}$. By \cref{finite}, $\mathscr O_{n+2}/I$ is always finitely generated over $\mathscr O_1$.%Recolocación del FT(pi,g)

%\begin{theorem}\label{samuel} Let $I=d\pi\left(\Der(-\log G)\right)$. Then,
%\[
%\mu_I(f)=e\left((s); \frac{\mathscr O_{n+2}}{FT(\pi,G)}\right),
%\]
%the Samuel multiplicity of $\mathscr O_{n+2}/I$ with respect to the maximal ideal $(s)$.
%\end{theorem}
\begin{theorem}\label{samuel} The image Milnor number of $f$ equals the Samuel multiplicity of $\mathscr O_{n+2}/FT(\pi,G)$ with respect to the maximal ideal $(s)$, i.e.,
\[
\mu_I(f)=e\left((s); \frac{\mathscr O_{n+2}}{FT(\pi,G)}\right).
\]
\end{theorem}

\begin{proof} Take $s_0\ne0$ close enough to the origin in $\CC$. By the conservation of the multiplicity (see, for example, \cite[Corollary E.5]{Mond2020}),
\begin{align*}
e\left((s); \frac{\mathscr O_{n+2}}{FT(\pi,G)}\right)&=\sum_{(y,s_0)\in V(FT(\pi,G))} e\left((s-s_0); \frac{\mathscr O_{n+2,(y,s_0)}}{FT(\pi,G)}\right)\\
&=\sum_{(y,s_0)\in V(FT(\pi,G))} e\left((s-s_0); \frac{\mathscr O_{n+2,(y,s_0)}}{\left(\frac{\partial G}{\partial y_1},\dots,\frac{\partial G}{\partial y_{n+1}}\right)}\right)\\
&=\sum_{y\notin X_{s_0}} \dim_\CC \frac{\mathscr O_{n+1,y}}{\left(\frac{\partial g_{s_0}}{\partial y_1},\dots,\frac{\partial g_{s_0}}{\partial y_{n+1}}\right)}\\
&=\sum_{y\notin X_{s_0}} \mu(g_{s_0};y)\\
&=\mu_I(f).
\end{align*}
The second equality follows from \cref{V(I)}, the third one holds because $$\frac{\mathscr O_{n+2,(y,s_0)}}{\left(\frac{\partial G}{\partial y_1},\dots,\frac{\partial G}{\partial y_{n+1}}\right)}$$ is Cohen-Macaulay, and the last one is a consequence of a theorem due to Siersma, \cite[Theorem 2.3]{Siersma1991}.
\end{proof}

\begin{remark}\label{intersection} The multiplicity $e\left((s); \frac{\mathscr O_{n+2}}{FT(\pi,G)}\right)$ can be interpreted geometrically as a local intersection number
\[
e\left((s); \frac{\mathscr O_{n+2}}{FT(\pi,G)}\right)=i\left(j,V\big(FT(\pi,G)\big)\right),
\]
where $j\colon(\CC^{n+1},0)\to(\CC^{n+2},0)$ is the embedding $y\mapsto(y,0)$. We refer to Fulton's book \cite{Fulton1998} for details about the connection between the algebraic multiplicity and the local intersection number. 

An alternative explanation of this equality can be given as we know that the points of $V\big(FT(\pi,G)\big)$ are precisely the instability of $f_0$ and the critical points of $g_s$, where $g_s$ is the equation of $\im(f_s)$, that are not contained in the image of $f_s$, by \cref{V(I)}. Moreover, we have already mentioned that a result of Siersma (see \cite[Theorem 2.3]{Siersma1991}) says that the sum of the Milnor numbers of these critical points is equal to the image Milnor number of $f_0$. In other words,
$$\sum_{y\notin X_{s_0}} \mu(g_{s_0};y)=\mu_I(f), $$
which is equal to the intersection multiplicity $i\big(j,V\big(FT(\pi,G)\big)\big)$ and, in turn, equal to the Samuel multiplicity $e\big((s); \frac{\mathscr O_{n+2}}{FT(\pi,G)}\big)$ (by \cref{samuel} above). See \cref{rem:siersma+curva,fig:SiersmaesDios} below.
\end{remark}

\begin{corollary} In terms of \cref{samuel},
\[
\mu_I(f)\le \dim_\CC\frac{\mathscr O_{n+2}}{FT(\pi,G)+(s)},
\]
with equality if and only if $\mathscr O_{n+2}/FT(\pi,G)$ is Cohen-Macaulay of dimension $1$.%Está bien, no falta (s). Esto es porque el modulo que pones en la multiplicidad de Samuel e(*,M) es ese mismo, que si es C-M entonces se cumple la conservación de multiplicidad mirando esas dimensiones
\end{corollary}

The following definition is an adaptation of the definition of the logarithmic characteristic variety $LC(\mathcal X)$ introduced by Saito in \cite{Saito1980}, where we consider the module $\Der(-\log G)$ instead of $\Der(-\log \mathcal X)$. 

Let $T^*\CC^{n+2}$ be the cotangent bundle of $\CC^{n+2}$. Given an open set $U\subset \CC^{n+2}$, $T_U^*\CC^{n+2}$ is the restriction of $T^*\CC^{n+2}$ to $U$. An element of $T_U^*\CC^{n+2}$ will be of the form $(y,s;\alpha)$, where $(y,s)\in U$ and $\alpha\colon\CC^{n+2}\to\CC$ is a linear form.

Given a holomorphic function $h\colon U\to\CC$ (or a germ $h\colon(\CC^{n+2},0)\to\CC$), we denote by $Dh\colon U\to T^*\CC^{n+2}$ (resp. $Dh\colon(\CC^{n+2},0)\to T^*\CC^{n+2}$) the differential, that is, the section of $T^*\CC^{n+2}$ given by $Dh(y,s)=\big(y,s;dh_{(y,s)}\big)$.

\begin{definition} Assume that $\xi_1,\dots,\xi_r$ generate $\Der(-\log G)$ on some open neighbourhood $U$ of the origin in $\CC^{n+2}$. Then, the \emph{logarithmic characteristic variety} of $G$ is defined as follows:%Es finitamente generado porque es submodul ode un f.g. y es Noetheriano
\[
LC_U(G)=\left\{(y,s;\alpha)\in T_U^*\CC^{n+2}:\ \alpha(\xi_i(y,s))=0,\ \forall i=1,\dots,r\right\}.
\]
The variety $LC(G)$ is the germ of $LC_U(G)$ along $T_0^*\CC^{n+2}$.
\end{definition}

It is easy to see that
\begin{equation}\label{eq:Dpi-1}
D\pi^{-1}\big(LC(G)\big)=V\big(FT(\pi,G)\big).
\end{equation}

To give the equations of $LC(G)$, suppose that the module $\Der(-\log G)$ is generated by germs of vector fields $\xi_1,\dots,\xi_r$ and that
\[
\xi_i=\sum_{i=1}^{n+1}a_j^i \frac{\partial}{\partial y_j}+b_i \frac{\partial}{\partial s},
\]
for some $a_j^i,b_i\in\mathscr O_{n+2}$. Denote by $(y_1,\dots,y_{n+1},s,p_1,\dots,p_{n+1},q)$ the coordinates of $T^*\CC^{n+2}$.
Then $LC(G)$ has equations $\xi^*_i=0$, $i=1,\dots,r$, where
\[
\xi_i^*=\sum_{i=1}^{n+1}a_j^i p_j+b_i q\in \mathscr O_{n+2}[p_1,\dots,p_{n+1},q].
\]
This shows that $LC(G)$ is independent of the choice of the neighbourhood $U$. It is not difficult to see that it is also independent of the choice of the generators. We remark that $LC(G)$ is considered with the possibly non-reduced structure given by the ideal generated by $\xi^*_i$, $i=1,\dots,r$.

We know that $LC(\mathcal X)$ has dimension $n+2$, for $(\mathcal X,0)$ is holonomic  (see \cite[Proposition 1.14]{Bruce1988}). We compute the dimension of $LC(G)$ in the next proposition.

\begin{proposition}\label{dimLC(G)} The variety $LC(G)$ has dimension $n+3$. Furthermore, the variety $V\big(FT(\pi,G)\big)$ has dimension $1$ if, and only if, $f$ is not stable and is empty otherwise.
\end{proposition}

\begin{proof} By \cref{finite}, 
\[
\dim_\CC\frac{\mathscr O_{n+2}}{FT(\pi,G)+(s)}<\infty,
\]
so $\codim V\big(FT(\pi,G)\big)\geq n+1$, and, by \cref{eq:Dpi-1}, $\codim LC(G)\geq n+1$. Hence, $\dim V\big(FT(\pi,G)\big)\leq1$ and $\dim LC(G)\leq n+3$.%La dimensión es esa porque como el cociente es C-ev finito contiene a una potencia de m, por lo tanto ese conjunto es solo el 0, pero eso es una hipersuperficie del V(FT), por lo tanto este último tiene dimensión 1 a lo sumo.

On the other hand, $G$ has isolated critical value at the origin by the curve selection Lemma. We take a small enough open neighbourhood $U$ of the origin in $\CC^{n+2}$ such that $0$ is the only critical value of $G$ on $U$. For each $(y,s)\in U- \mathcal X$, $G$ is regular at $(y,s)$. We use the argument given in the proof of \cref{V(I)} and deduce that $\Der(-\log G)$ is generated in a neighbourhood of $(y,s)$ by $n+1$ vector fields $\xi_1,\dots,\xi_{n+1}$. It follows that $LC(G)$ is given in a neighbourhood of $(y,s)$ by $n+1$ equations $\xi^*_i=0$, $i=1,\dots,n+1$, so $\codim LC(G)\le n+1$. Hence, $\dim LC(G)=n+3$ and $\dim V\big(FT(\pi,G)\big)=1$ whenever it is not empty. However, by \cref{finite},  $FT(\pi,G)$ contains the constants if, and only if, $f$ is stable, and $V\big(FT(\pi,G)\big)$ would be empty.
\end{proof}

Now we prove a weak version of Mond's conjecture \cref{MondConjecture} for the case $\mu_I(f)=0$. This is a generalization of \cite[Theorem 3.9]{GimenezConejero2022}, which is stated for corank one germs.

\begin{theorem}\label{WMC} We have that $\mu_I(f)=0$ if and only if $f$ is stable.
\end{theorem}

\begin{proof} Assume that $f$ is not stable. We know from \cref{dimLC(G)} that $\dim LC(G)=n+3$ and $\dim V\big(FT(\pi,G)\big)= 1$.
In other words, $$\dim \mathscr O_{n+2}/FT(\pi,G)=1=\dim\mathscr O_1$$ and, hence,
\[
\mu_I(f)=e\left((s);\frac{\mathscr O_{n+2}}{FT(\pi,G)}\right)\ge 1,
\] 
(see, for instance, \cite[Formula 14.2]{Matsumura1986}).
\end{proof}

\begin{remark}\label{rem:siersma+curva}
There is an equivalent way of proving this result using the geometry of $V\big(FT(\pi,G)\big)$.
On one hand, we know that
$$\mu_I(f)=i\left(j,V\big(FT(\pi,G)\big)\right),$$
as shown in \cref{intersection}. On the other, this intersection number is positive, because $V\big(FT(\pi,G)\big)$ is curve (see \cref{dimLC(G)}), so the result follows. See \cref{fig:SiersmaesDios} for a general overview of this reasoning.
\end{remark}

\begin{figure}
	\centering
		\includegraphics[width=1.00\textwidth]{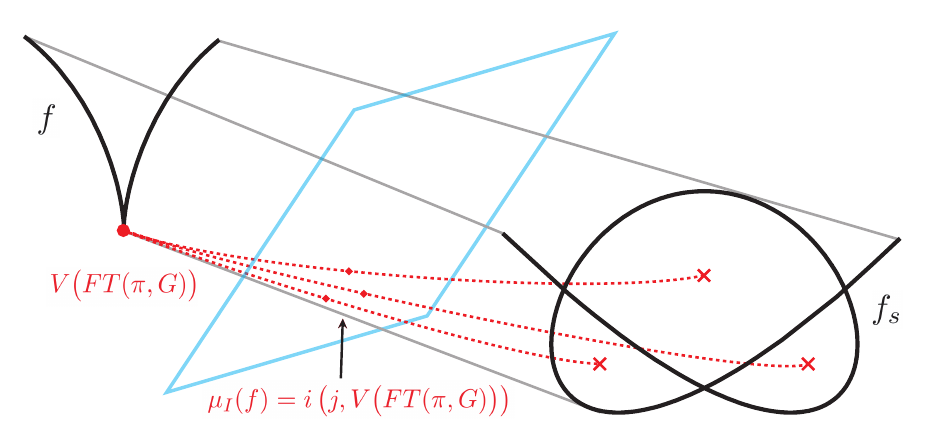}
	\caption{The curve $V\big(FT(\pi,G)\big)$ and the intersection number $i\left(j,V\big(FT(\pi,G)\big)\right)$.}
	\label{fig:SiersmaesDios}
\end{figure}

\begin{corollary}\label{MCae1}
Mond's conjecture is true for germs of $\eqA_e$-codimension 1 (see \cref{MondConjecture}).
\end{corollary}
\begin{proof}
The general case follows from \cref{WMC}. If $f$ is quasi-homogeneus then its $\eqA_e$-codimension coincides with the dimension as vector space of the Jacobian module $M(g)$, see \cite{Bobadilla2019}, hence
$$ 0<\mu_I(f)\leq \dim_\C M(g)=1,$$
by \cite[Proof of Theorem 6.1]{Bobadilla2019}.
\end{proof}

Cooper, Mond and Wik Atique proved this result for corank one germs in \cite[Theorem 7.2]{Cooper2002}. Actually, they proved that, in corank one, germs with $\eqA_e$-codimension 1 have image Milnor number equal to 1. This fact led Houston to prove that Mond's conjecture holds for a class of germs called \textit{augmentations} of corank one germs that have $\eqA_e$-codimension one (see \cite[Corollary 6.8]{Houston2002}). We can generalize the proof to germs of any corank using the same ideas. The concept of augmentation was first introduced in \cite{Goryunov1983}, and the reader can find its definition in, for example, \cite[Definition 3.1]{Houston1998}.

\begin{corollary}
Suppose that $A_{F,g}(f)$ is an augmentation with $g$ of a germ $f:\GS{n}{n+1}$ of $\eqA_e$-codimension one. If $f$ or $g$ are quasi-homogeneous, $A_{F,g}(f)$ satisfies Mond's conjecture and, more precisely,
$$\aecodim\big(A_{F,g}(f)\big)=\aecodim(f)\tau(g)\leq \mu_I(f)\mu(g)=\mu_I\big(A_{F,g}(f)\big),$$
where $\tau$ denotes the Tjurina number, with equality if $f$ and $g$ are quasi-homogeneous.
\end{corollary}
\begin{proof}
The result is a consequence of \cref{MCae1} and \cite[Theorem 6.7]{Houston2002} (see also \cite[Theorem 3.3]{Houston1998}, which controls the $\eqA_e$-codimension of augmentations). There is, however, a small consideration to be made: the equality 
$$ \mu(g)\mu_I(f)=\mu_I\big(A_{F,g}(f)\big) $$
is stated for corank one monogerms in \cite[Corollary 6.4]{Houston2002}. One can give the same proof for multigerms of any corank, using the well known fact that the image of a stable perturbation of any $f$ with our hypothesis has the homotopy type of a wedge of spheres (see, for example, \cite[Proposition 8.3]{Mond2020}).
\end{proof}

\begin{theorem}\label{thm:muinterseccion} Let $f:(\CC^n,S)\to(\CC^{n+1})$ be a germ such that $(n,n+1)$ are nice dimensions or it has corank one. Then,
\[
\mu_I(f)=i\big(D\pi\circ j, LC(G)\big),
\]
where $j(y)=(y,0)$, as in \cref{intersection}.
\end{theorem}
\begin{proof}
In fact, we have
\[
i\big(D\pi\circ j,LC(G)\big)=e\left(\mathfrak m_{n+3,(0,\pi)};\mathscr O_{LC(G),D\pi(0)}\right),
\]
where we consider $\mathscr O_{LC(G),D\pi(0)}$ as a module over $\mathscr O_{n+3,(0,\pi)}$ via the projection $T^*\CC^{n+2}\to\CC\times\CC^{n+2}=\CC^{n+3}$ given by $(y,s;\alpha)\mapsto(s;\alpha)$.

As in the proof of \cref{samuel} we use the conservation of the multiplicity (see, for example, \cite[Corollary E.5]{Mond2020}). For $(s_0;\pi)$ with $s_0\neq 0$,
\[
e\left(\mathfrak m_{n+3,(0,d\pi)};\mathscr O_{LC(G),D\pi(0)}\right)=
\sum_{D\pi(y,s_0)\in LC(G)} e\left(\mathfrak m_{n+3,(s_0,\pi)};\mathscr O_{LC(G),D\pi(y,s_0)}\right).
\]
Since $G$ is regular at each $(y,s_0)$, $\Der(-\log G)$ is generated by $n+1$ vector fields at $(y,s_0)$, therefore $LC(G)$ is given by $n+1$ equations at $D\pi(y,s_0)$. So, $LC(G)$ is Cohen-Macaulay at $D\pi(y,s_0)$ and, hence,
\begin{align*}
	e\left(\mathfrak m_{n+3,(s_0,\pi)};\mathscr O_{LC(G),D\pi(y,s_0)}\right)&=\dim_\C\frac{\mathscr O_{LC(G),D\pi(y,s_0)}}{\mathfrak m_{n+3,(s_0,\pi)}\mathscr O_{LC(G),D\pi(y,s_0)}}\\
	&=\dim_\C\frac{\mathscr O_{V(FT(\pi,G)),(y,s_0)}}{\mathfrak m_{1,s_0}\mathscr O_{V(FT(\pi,G)),(y,s_0)}}\\
	&=\dim_\C\frac{\mathscr O_{n+2,(y,s_0)}}{FT(\pi,G)+\mathfrak m_{1,s_0}}\\
	&= \mu(g_{s_0};y);
\end{align*}
where the last equality is given by \cref{V(I)}.

Using \cite[Theorem 2.3]{Siersma1991} as in \cref{samuel}, this shows that $$e\left(\mathfrak m_{n+3,(0,\pi)};\mathscr O_{LC(G),{D\pi}(0)}\right)=\mu_I(f),$$ and the result follows.
\end{proof}

Observe that, as we are using multiplicities in the proof of \cref{thm:muinterseccion}, the equality
\[
\mu_I(f)=\dim_\CC\frac{\mathscr O_{n+2}}{FT(\pi,G)+(s)},
\]
holds if and only if $LC(G)$ is Cohen-Macaulay at $(0,\pi)$.
This leads to the following result.

\begin{theorem}
Mond's conjecture is true for germs $f$ with one parameter stable unfoldings $F$ (\textsc{OPSU}) and such that $LC(G)$ is Cohen-Macaulay, for $G$ given by $F$ (see \cref{MondConjecture}).
\end{theorem}
\begin{proof}
As we were saying, we have conservation of the multiplicity, since $LC(G)$ is Cohen-Macaulay, so
\begin{equation}\label{eq:muidim}
\mu_I(f)=\dim_\CC\frac{\mathscr O_{n+1}}{FT(\pi,G)+(s)}.
\end{equation}

Now, by \cite[Theorem 8.7]{Mond2020} (see the original version in \cite{Damon1991a}),
$$T^1\eqA_e(f)\cong\frac{\theta(i)}{ti(\theta_{n+1})+i^*\big(\Der(-\log \mathcal{X})\big)},$$
where $i$ is given by the commutative diagram
$$\begin{tikzcd}
(\CC^n\times\CC,S\times\left\{0\right\})\arrow[r,"F"]&(\CC^{n+1}\times\CC,0)\\
(\CC^n,S)\arrow[r,"f"]\arrow[u,"j",hook]&(\CC^{n+1},0)\arrow[u,"i"',hook]
\end{tikzcd},$$
and $i$ is transverse to $F$.
However, we can do the identifications $\theta(i)/di(\theta_Y)\cong\theta(\pi)/\mathfrak{m}_1\theta(\pi)$ and $\theta(\pi)\cong\mathscr{O}_{n+1}$, so 
\begin{equation}\label{eq:aedim}\begin{aligned}T^1\eqA_e(f)&\cong \frac{\theta(\pi)}{d\pi\big(\Der(-\log \mathcal{X})\big)+\mathfrak{m}_1\theta(\pi)}\\
&\cong \frac{\mathscr{O}_{n+1}}{d\pi\big(\Der(-\log \mathcal{X})\big)+(s)}.\end{aligned}\end{equation}

Recall that, by definition, $FT(G,\pi)=d\pi \big(\Der(-\log G)\big)$, and the decomposition 
$$d\pi\big(\Der(-\log\mathcal X)\big)=d\pi\big(\Der(-\log G)\big)+(s)$$ in the quasi-homogeneous case given in \cref{euler2}. Hence, comparing \cref{eq:muidim,eq:aedim}, the result follows.
\end{proof}

\section{Bifurcation set}\label{sec: bif}

In this section we prove that the bifurcation set of an $\eqA$-finite germ $f:\GS{n}{n+1}$ is a hypersurface that has pure dimension.  This is, indeed, the only case where it is not known whether the bifurcation set is a hypersurface (without any additional hypotheses). 

For germs $f:\GS{n}{p}$ with $n\geq p$, this property is shown in \cite[Theorem 8.8]{Mond2020}, whose proof relies on the fact that the discriminant is a free divisor in this case (see \cite[Corollary 6.13]{Looijenga1984}, cf. \cite[Proposition 8.9]{Mond2020}). This is no longer true for $p>n$.%certain object being Cohen-Macaulay\todo{También estaba lo de los liftables eran un modulo libre, pero creo que esto también es una razón, no?} (which is only known in these dimensions). 
For germs with $p\geq n+2$ the bifurcation set can have greater codimension, see \cite[Example 9.5]{Mond2020}. Finally, when $p=n+1$, it was only known that the bifurcation set is a hypersurface for germs of corank one (see \cite[Proposition 9.15]{Mond2020}). This, in turn, relied in the good structure of the \textit{multiple point spaces}, which lose this good behaviour when we study germs with greater corank.
\newline

\begin{figure}[ht]
	\centering
		\includegraphics[width=0.85\textwidth]{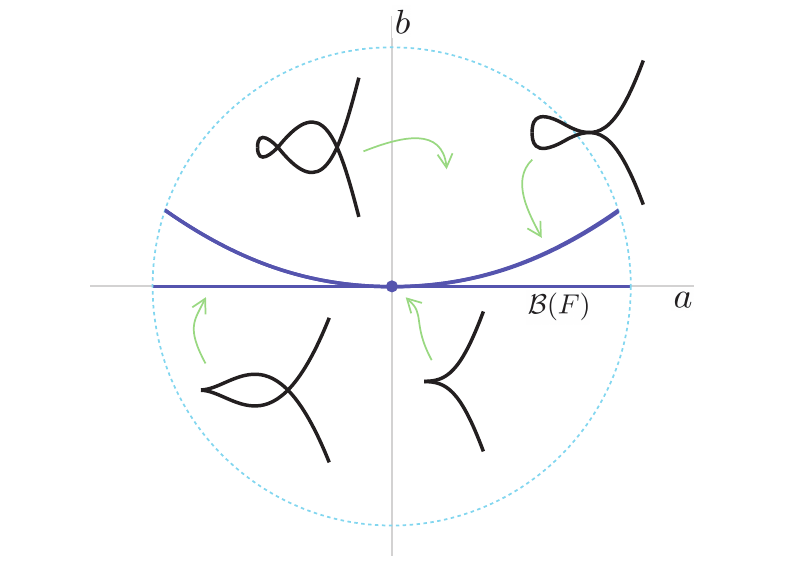}
	\caption{Representation of the parameter space of the miniversal unfolding $F(x,a,b)=(x^2,x^5+ax^3+bx,a,b)$ of $f(x)=(x^2,x^5)$.}
	\label{fig: versal x2x5}
\end{figure}

Consider a map germ $f:(\CC^n,S)\to(\CC^{n+1},0)$ with $\eqA_e$-codimension $d>0$ and a miniversal unfolding (see \cref{fig: versal x2x5})
\begin{align*}
	F:\big(\CC^n\times\CC^d,S\times\{0\}\big)&\to(\CC^{n+1}\times\CC^d,0)\\
	(x,\lambda)&\mapsto \big(f_\lambda(x),\lambda\big).
\end{align*}

We want to use some ideas Goryunov used in \cite{Goryunov1991} to reach a version for more dimensions of his \cite[Corollary]{Goryunov1991} given for germs from $\CC^2$ to $\CC^3$. Indeed, it is implicit in Goryunov's work the result given here as \cref{WMC} or, as it is proven in those dimensions, Mond's conjecture (see \cite[Theorem 4.2]{deJong1991} and \cite[Theorem 2.3]{Mond1995}). 

We can consider a generic line $\tilde{\ell}$ in $\CC^d$ through the origin to compute the multiplicity of the bifurcation set $\mathscr{B}(F)$. Hence, when we shift the line $\tilde{\ell}$ to the new line $\ell$, the number $\nu$ of points in $\ell\cap\mathscr{B}(F)$ is the multiplicity of $\mathscr{B}(F)$ at the origin. Observe, however, that this multiplicity is positive if, and only if, $\mathscr{B}(F)$ is a hypersurface and, to prove that it is positive, we will use \cref{WMC}. 

%To fix notation, say that the intersection between $\ell$ and $\mathscr{B}(F)$ is the possibly empty set $\left\{p_1,\dots,p_r\right\}$.

We can consider the pullback $F^\ell$ induced by $\ell$ and $F$, i.e., $F^\ell(x,t)\coloneqq\big(f_{\ell(t)}(x),t\big)$ where $\ell(t)$ parametrizes the line $\ell$. Notice that this is a deformation of the pullback $F^{\tilde{\ell}}$ induced by $\tilde{\ell}$ and $F$.

\begin{theorem}\label{thm: hypersurface}
If $f$ has corank one or it is in the nice dimensions, the bifurcation set $\mathscr{B}(F)$ is a hypersurface.
\end{theorem}
\begin{proof}
We are going to use stratified Morse theory with a function induced from the projection $\pi$ to the parameter $\lambda$ in the image of $F^\ell$:
$$
\begin{tikzcd}[column sep=10pt,row sep=1pt]
\phi:\im F^\ell \arrow[r,]& \ell\\
\big(f_{\ell(t)}(x),t\big)\arrow[r,mapsto]& \left|\ell(t)-q_0\right|
\end{tikzcd},
$$
where $q_0$ is a generic point of $\ell$, which can be considered to lie in $\ell-\mathscr{B}(F)$. 
We also use the stratification given by the stable types (and the isolated unstable points) of $F^\ell$. %For the sake of simplicity, let us assume that $\phi$ is a Morse function, otherwise the argument is very similar, but taking a morsification of $\phi$.

Observe that a critical point of $\pi$ induces a critical point of $\phi$. Furthermore, the points $p_i$ of the intersection $\ell\cap\mathscr{B}(F)$ are, precisely, the critical points of $\pi$ by \cref{isolated} (or, to be more precise, a very easy adaptation of this result to $F^\ell$). By contradiction, if the intersection $\ell\cap\mathscr{B}(F)$ were empty, there would be no critical points, and, by Morse theory, the image of $F^\ell$ would be a deformation retract of the image of $f_{q_0}$. But this is absurd. Indeed, $F^\ell$ is a deformation of $F^{\tilde{\ell}}$, and $f_{q_0}$ is a stable perturbation of $f$ (so it has the homotopy type of a wedge of $\mu_I(f)$ spheres of dimension $n$). As $f$ has an instability, the number $\mu_I(f)$ is positive (by \cref{WMC}), but the image of $F^\ell$ has only non-trivial homology in dimension zero and, possibly, in dimension $n+1$.
\end{proof}

With these ideas, we can refine the result.

\begin{theorem}
In the conditions of \cref{thm: hypersurface},
$$ m\big(\mathscr{B}(F)\big)\leq \mu_I(f)+\mu_I\big(F^{\tilde{\ell}}\big). $$
\end{theorem}
\begin{proof}
Using the exact sequence of the pair,
$$ 0\to H_{n+1}\big( \im F^\ell\big)\to H_{n+1}\big( \im F^\ell,\im f_{q_0}\big)\to H_n(f_{q_0})\to0,$$
we only need to confirm that each point $p_i$ (i.e., each critical point of $\pi$) contributes with at least one copy of $\ZZ$ to $H_{n+1}\big( \im F^\ell,\im f_{q_0}\big)$. Indeed, if $\phi$ is not a stratified Morse function we can consider a Morsification, which has at least one critical point for each critical point of $\phi$ (hence $\pi$). As these images are analytic and $\phi$ is the module of a complex analytic function, the tangential Morse data is %esto está en el libro, mirar 267 (módulo de una compleja, mirar diferencial y ver que iv es autovector de -1)
$$ (D^m,\partial D^m), $$
where $m$ is the dimension of the stratum that contains the critical point. Furthermore, as we deal with hypersurfaces, a theorem of Lê (see \cite{Trang1979}, but also the more easy to access \cite[pp. 187-188]{Goresky1983}) says that the normal Morse data is homotopic to
$$ \left(\vee D^{n+1-m},\partial \vee D^{n+1-m}\right).$$
Hence, the Morse data is
\begin{align*}
	(D^m,\partial D^m)&\times\left(\vee D^{n+1-m};\partial \vee D^{n+1-m}\right)\\
	&=\left(D^m\times\vee D^{n+1-m};\partial D^m\times\vee D^{n+1-m}\cup  D^m\times\partial\vee D^{n+1-m} \right)\\
	&\simeq \left(\vee D^{n+1};C_{<n+1}\right),
\end{align*}
where $C_{< n+1}$ is a $CW$-complex of dimension lower than $n+1$. The result follows from here, counting the new cells per critical point.
\end{proof}

\begin{proposition}
In the conditions of \cref{thm: hypersurface}, the bifurcation set is pure dimensional.
\end{proposition}
\begin{proof}
By contradiction, assume that $\mathscr{B}(F)$, with reduced structure, has a component that is not a hypersurface. We can consider a point $\lambda$ that lies in that component but not in the other components. By the openness of versality (see \cite[Theorem 3.7]{Wall1981}, cf. \cite[Theorem 5.6]{Mond2020}), $F$ is also a versal unfolding of $f_\lambda$ at each point. This is already a contradiction, as it is not a hypersurface.
\end{proof}

\begin{remark}
It is important to notice that the bifurcation set considered with the reduced structure is not, in general, irreducible. Indeed, in most cases it will have more than one irreducible component (see, for example, \cite[Example 5.8]{Mond2020}).
\end{remark}

\bibliographystyle{plain}
\bibliography{FullBib}

\end{document}